\documentclass[12pt,twoside]{article}
\usepackage{amsmath}
\usepackage{amsthm}
\oddsidemargin=0cm \evensidemargin=0mm \newfont{\titre}{cmbx12 at 16pt}

\newtheorem{lem}{Lemma}[section]
\newtheorem{thm}{Theorem}[section]
\newtheorem{prop}{Proposition}[section]

\small

\def\mathr{{\sf I\kern-0.1em R}}
\def\Mathz{{\Large \kern0.3em\sf \Large Z\kern-0.4em Z\kern0.3em}}
\def\mathz{{\kern0.3em\sf Z\kern-0.8em Z\kern0.2em}}
\def\mathc{{\kern0.1em\sf C\kern-0.8em C\kern0.2em}}

\def\11{I\!\!I}  \def\RR{I\!\!R} \def\NN{I\!\! N}
\def\RR{{\sf I\kern-0.1em R}}
\title{Feedback Control of Nonlinear Dissipative Systems
by  Finite Determining   Parameters \\
- A Reaction-diffusion Paradigm}
\author{{Abderrahim ~Azouani$^1$ and Edriss S. ~Titi} $^{2, 3}$}
\date{\small  $^{1}$Freie Universit\"at Berlin, Institut f\"ur Mathematik I,
Arnimallee 7, Berlin, Germany.\\
{\em azouani@math.fu-berlin.de}\\
$^{2}$ Departments of Mathematics and
of Mechanical and Aerospace Engineering, University of California, Irvine,
CA 92697-3875, USA.\\
$^{3}$  Department of Computer Science and Applied Mathematics,
 Weizmann Institute of Science, Rehovot 76100, Israel.\\
{\em etiti@math.uci.edu} and {\em edriss.titi@weizmann.ac.il}
\\
May 22, 2014}

\begin{document}
\maketitle{}

\abstract{

We introduce here a simple finite-dimensional feedback control scheme for stabilizing solutions of infinite-dimensional dissipative evolution equations, such as reaction-diffusion systems, the Navier-Stokes equations and the Kuramoto-Sivashinsky equation. The designed feedback control scheme takes advantage of the fact that such systems possess finite number of determining parameters (degrees of freedom), namely,  finite number of  determining Fourier modes, determining nodes, and determining interpolants and projections. In particular, the  feedback control scheme uses finitely many of such observables and controllers. This observation is of a particular  interest since it implies that our approach has far more reaching applications, in particular, in data assimilation. Moreover, we emphasize that our scheme treats all kinds of the determining projections, as well as, the various dissipative equations  with one unified approach. However, for  the sake of simplicity we demonstrate  our approach in this paper to  a one-dimensional reaction-diffusion equation paradigm. }
\\
\\
{\bf Keywords.} Reaction-diffusion, Navier-Stokes equations, feedback control, data assimilation, determining modes, determining
 nodes, determining volume elements.\\
{\bf Mathematics Subject Classification (2000):} 35K57, 37L25, 37L30, 37N35, 93B52, 93C20, 93D15.

\section{Introduction}
Dissipative dynamical systems, such as the Navier-Stokes equations, the
Kuramoto-Sivashinsky equation, the complex Ginzburg-Landau equation and various reaction-diffusion systems are known to have a finite-dimensional asymptotic (in time) behavior  (see, e.g.,  \cite{CJT,Cock97, Cos:nse, FMRT, Hale88, Robinson, Sell-You, Tem97},
 and references therein). This is evident due to the fact that such systems possess  finite-dimensional global attractors (\cite{BabVish92, Cos:nse, Cos88, Robinson,Sell-You,Tem97}), and finite number  of determining modes (\cite{Foias-Prodi, Foias-Manley-Temam-Treve, FMRT, Jones93}), determining nodes (\cite{FMRT, Foias-Temam-nodes1, Foias-Temam-nodes2, Foias91, Jo:nse2D, Jones93, Kukavica}), determining volume
elements (\cite{Foias91,Jones92B}) and other finite number of determining parameters (degrees of freedom)  such as finite elements and other interpolation polynomials (\cite{CJT,Cock97,Foias-Temam-nodes2}.)
Moreover, some of these systems,  which enjoy the property of separation of spatial scales,  are also known to have a
finite dimensional inertial manifolds (see, e.g., \cite{Cos:nse, Cos88, Foias88, Foias-Sell-Titi, Tem97}, and references therein). That is, in the presence of separation of spatial scales the long-term dynamics of such a system is equivalent to that of a finite system of ordinary differential equations.

There has been some interesting work on reduction methods, with applications focused on  scientific computing and feedback control theory, taking advantage of the finite-dimensional asymptotic behavior of these dissipative dynamical systems  (see, e.g., \cite{Christofides2000,Christofides2003,Foias-Jolly-Kevrekidis-Sell-Titi,Jolly-Kevrekidis-Titi, Yannis-Control} and references therein) . However, there has been very little  rigorous  analytical work, in particular in  the context of feedback control theory, justifying these applications.
In the case of separation of spatial scales,  and hence the existence of inertial manifolds, the authors of \cite{Ros:rd} and \cite{Ros:rd-im}  provide an example  of finite-dimensional feedback control (lumped feedback control) that drives  the dynamics of  one-dimensional reaction-diffusion system to an a priori specified finite-dimensional dynamics. It is worth stressing again that in the case of inertial manifold the dynamics of the underlying evolution equation is  equivalent to that of an ordinary differential  equations to begin with. However, the main challenge is in being able to  provide a representation of this ODE system in the relevant parameters dictated by the applications. In \cite{Foias91} and \cite{Cock97} the authors have shown that if a certain dissipative system has separation of scales, and hence an inertial manifold, then such a manifold can be parameterized by any set of adequate parameters, e.g.  Fourier modes, nodal values, local volume averages, etc... In the above mentioned work of \cite{Ros:rd} and \cite{Ros:rd-im} the authors employed such an equivalence in the parameterization of the inertial manifolds to show their results.

In this paper we  propose a new feedback control for controlling general dissipative evolution equations using any of the determining systems of parameters (modes, nodes, volume elements, etc...) without requiring the presence of separation in spatial scales, i.e. without assuming the existence of an inertial manifold. To fix ideas we demonstrate our idea for a simple reaction diffusion equation, the Chafee-Infante equation, which is the real Ginzburg-Landau equation. It is worth mentioning, however,  that this new idea has a far more reaching areas of applications, other than feedback control, such as  in data assimilations for weather prediction \cite{Azouani-Olson-Titi,Bessaih-Olson-Titi}.  In addition, one can use this approach to show that the long time-time dynamics of the underlying dissipative evolution equation, such as the two-dimensional Navier-Stokes equations, can be imbedded in an infinite-dimensional dynamical system that is induced by  an ordinary differential equations, named {\it determining form}, which is governed by a globally Lipschitz vector field,  cf. \cite{Foias-Jolly-Kravchenko,Foias-Jolly-Kravchenko-Titi1} and \cite{Foias-Jolly-Kravchenko-Titi2}.

In this paper we will use the Chafee-Infante reaction-diffusion equation
\begin{equation}
\frac{\partial u}{ \partial t}-\nu \, u_{xx}-\alpha u+u^3=0 \label{RD}
\end{equation}
\begin{equation}
u_x (0)=u_x (L)=0
\label{bc1}
\end{equation}
for $\alpha > 0$, large enough, as a paradigm to fix ideas and to
use the notions of finite number of
determining modes, nodes and volume elements to design feedback control
to stabilize the $\mathbf{v} (x) \equiv 0$ unstable steady state solution
of (\ref{RD})-(\ref{bc1}).
Indeed, by linearizing
equation (\ref{RD}) about $\mathbf{v} \equiv 0$ one obtains the
linear equation
\begin{eqnarray}
\frac{\partial \mathbf{v}}{ \partial t}-\nu \mathbf{v}_{xx}-
\alpha \mathbf{v}=0 \label{RD-st} \\
\mathbf{v}_x (0)=\mathbf{v}_x (L)=0 \nonumber
\end{eqnarray}
We solve equation (\ref{RD-st}) with initial condition
$\mathbf{v}_0(x)=A_k \, \cos(\frac{kx}{L} \pi)$, where $A_k \in \RR $,
by seeking  a solution of the form
$\mathbf{v}(x, t)=a_k(t) \, \cos(\frac{k \, x}{L} \pi)$.
Therefore, we obtain
\begin{equation*}
\dot{a_k}+ \nu a_k \, (\frac{\pi k} {L})^2 -\alpha \, a_k=0, \quad \hbox{with} \quad a_k(0)=A_k;
\end{equation*}
and whose solution is
\begin{equation*}
a_k(t)=A_k \, e^{\left(-\nu \, (\frac{\pi k} {L})^2+\alpha \right) \, t}.
\end{equation*}
Therefore, for $\alpha >0$, large enough, all the low wave numbers
$k^2 < \frac{\alpha \, L^2}{\pi^2 \nu}$ are unstable. Consequently,
the dimension of the unstable manifold of $\mathbf{v}\equiv0$ behaves
like $\sqrt{\frac{\alpha \, L^2} {\nu}}$ (see, for instance,
\cite{BabVish92, Hale88} and \cite{Tem97}
for a similar analysis).\\
The aim of this paper is to design a feedback control that stabilizes
$\mathbf{v} \equiv 0$, for example, either by observing the values of
the solutions at certain nodal points, local averages of the solutions
in subintervals of $[0, L]$, or by observing finitely many of their
Fourier modes. Based on the
above discussion, a naive analysis would suggest that one would need about
$\sqrt{\frac{L^2 \alpha}{\nu}}$ feedback controllers to stabilize
$\mathbf{v}\equiv 0$.\\
In this paper we will give a rigorous justification to this assertion.
First, we demonstrate our result for the case of local averages, which
is the most straightforward approach. Later, we present a more general
abstract result, that unifies our approach, utilizing all sorts of
approximate interpolant operators, as observables and controllers, and show that this
abstract approach applies to the Fourier modes, local volumes
(i.e. local averages) and nodal values as particular examples.
It is worth mentioning that the same feedback control scheme  can be used to stabilize any other time-dependent solution, $v(x,t)$, of (\ref{RD})-(\ref{bc1}). The details of the proof are similar to the ones presented here for  stabilizing the zero solution; thus,  for the sake of simplicity they will not be provided. Furthermore, similar scheme   can be also implemented for feedback control of other nonlinear dissipative dynamical systems, such as the two-dimensional Navier-Stokes equations, the Kuramoto-Sivashinsky equation and reaction-diffusion systems. A computational study concerning the implementation of this feedback control scheme for various nonlinear dissipative equations will be reported in a forthcoming work \cite{Lunasin-Titi}. In addition, one can design similar feedback control algorithm with stochastically noisy observables and controllers  to stabilize, in the average,  given solutions; within errors that are determined by the standard deviation  of the noise. This can be achieved by  combining  some of the  ideas presented in \cite{Bessaih-Olson-Titi} with those presented in the present paper, a subject of future work.

\section{Finite volume elements feedback control}
\label{sec:fin-vol-feed}
To fix ideas we propose the following feedback control system for
(\ref{RD})-(\ref{bc1}) in order
to stabilize the steady state solution $\mathbf{v}\equiv 0$,
\begin{eqnarray}
\frac{\partial u}{\partial t}- \nu \, u_{xx}-\alpha u+u^3 =-\mu
\sum_{k=1}^{N} \overline{u}_k \; \chi_{J_k\strut} (x) \label{av-RD} \\
u_x(0)=u_x(L)=0, \label{av-bc1}
\end{eqnarray}
where $J_k=\left[(k-1) \frac{L}{N}, k \frac{L}{N}\right),$
for $k=1,\dots,N-1$, and $J_N = \left[(N-1) \frac{L}{N},  L \right]$; moreover, $\chi_{J_k \strut}(x)$
is the characteristic function of the interval $J_k$, for $k=1,\dots,N$,  and
$$
\bar{\varphi}_k=\frac{1}{|J_k|} \, \int_{J_k} \, \varphi(x) \, \, dx=
\frac{N}{L} \, \int_{J_k} \, \varphi(x) \, \, dx.
$$
Here, the local averages of the solution, $\overline{u}_k$, for
$k=1,..., N,$ are the observables, and they are also used as the feedback
controllers
in (\ref{av-RD}). It is easy to observe $ \mathbf{v} \equiv 0$ is also a
steady state solution for
(\ref{av-RD})-(\ref{av-bc1}).\\
For $\varphi \in H^1([0, L])$ we define
\begin{equation}
\|\varphi\|_{H^1}^2:=\frac{1}{L^2} \, \int_{0}^{L}
\varphi^2(x)\, dx + \int_{0}^{L} \varphi^2_{x} (x) \, dx.
\label{norm-h1}
\end{equation}
Before showing that (\ref{av-RD})-(\ref{av-bc1}) globally stabilizes
the steady state $\mathbf{v}\equiv 0$,
one has to prove first the global existence and uniqueness of
the feedback system (\ref{av-RD})-(\ref{av-bc1}). In section
\ref{sec:exis-uni}, we will show in Theorem
\ref{thm:4.1} a result concerning global existence, uniqueness and stabilization for a general
family of finite-dimensional
feedback control that includes  system (\ref{av-RD})-(\ref{av-bc1}) as
a particular case.
Therefore, we will postpone this task of proving the global existence and
uniqueness until section
\ref{sec:exis-uni}, and we only
show here the global stability of $\mathbf{v} \equiv 0$. This is  in
order to fix ideas and to demonstrate our general approach.\\
Assuming the global existence and uniqueness of
(\ref{av-RD})-(\ref{av-bc1}), we will show in this section  that every solution
$u$ of (\ref{av-RD})-(\ref{av-bc1}) tends to zero, as
$t \rightarrow \infty$, under specific explicit assumptions on
$N, \nu, \alpha, L$ and $  \mu$
(see Theorem \ref{thm:fin-vol} for  details). But first we  need the following proposition to prove our result. We observe that similar
propositions were introduced and proved in \cite{Constantin96,Foias91, Jo:nse2D,Jones92B} and \cite{Jones93}
(see also \cite{Ros:rd} and \cite{Ros:rd-im}). We adapt here similar
ideas from \cite{Constantin96} for our proof.
\begin{prop}
\label{prop:89}
Let $\varphi \in H^1([0, L])$ then
\begin{equation}
\|\varphi(\cdot)-\sum_{k=1}^{N} \, \overline{\varphi}_k \, \chi_{J_k \strut}(\cdot)\|_{L^2} \,
\leq \, h \, \|\varphi_x\|_{L^2} \leq \, h \, \|\varphi\|_{H^1},
\label{8}
\end{equation}
where $h=\frac{L}{N}$. Moreover,
\begin{equation}
\|\varphi\|^2_{L^2}  \leq   h \,   \gamma^ 2(\varphi)+\left(\frac{h}{2 \, \pi}\right)^2 \, \|\varphi_{x}\|^2_{L^{2}},
\label{9}
\end{equation}
where $$\gamma^2(\varphi)=\sum_{k=1}^{N} \overline{\varphi}_k^2.$$
\end{prop}
{\bf{Proof.}}
\begin{eqnarray}
\|\varphi(\cdot)-\sum_{k=1}^{N} \, \overline{\varphi}_k \; \chi_{J_k \strut}(\cdot)\|^2_{L^2} \, &=& \int_{0}^{L} \left(\varphi(x)-\sum_{k=1}^{N} \,
\overline{\varphi}_k \, \chi_{J_k \strut}(x) \right)^2 dx \nonumber \\
&=& \int_{0}^{L}  \,\left(\varphi(x)  \sum_{k=1}^{N} \, \chi_{J_k \strut}(x)- \sum_{k=1}^{N} \, \overline{\varphi}_k \, \chi_{J_k \strut}(x)  \right)^2 dx,  \nonumber
\end{eqnarray}
where in the last equality we used the fact that $\sum_{k=1}^{N} \, \chi_{J_k \strut}(x) \equiv 1$.
Therefore,
\begin{eqnarray}
\|\varphi(\cdot)-\sum_{k=1}^{N} \, \overline{\varphi}_k \; \chi_{J_k \strut}(\cdot)\|^2_{L^2}
&=&\int_{0}^{L} \left(\sum_{k=1}^{N} \left(\varphi(x)-\overline{\varphi}_k \, \right) \,  \chi_{J_k \strut}(x))\right) \left(\sum_{l=1}^{N} \left(\varphi(x)-\overline{\varphi}_l \, \right) \, \chi_{J_l \strut}(x) \right) dx \nonumber \\
&=&\int_{0}^{L} \sum_{k,l=1}^{N} \, (\varphi(x)-\overline{\varphi}_k) \, (\varphi(x)-\overline{\varphi}_l) \, \chi_{J_k \strut}(x) \, \chi_{J_l \strut}(x) dx. \nonumber
\end{eqnarray}
Since $\chi_{J_l \strut}(x) \, \chi_{J_k \strut}(x) \equiv \chi_{J_k \strut}(x) \delta_{kl}$, it follows from the above that
\begin{eqnarray}
\|\varphi(\cdot)-\sum_{k=1}^{N} \, \overline{\varphi}_{k} \chi_{k \strut}(\cdot )\|_{L^2}^2&=&\int_{0}^{L} \sum_{k=1}^{N} \, \left(\varphi(x)-
\overline{\varphi}_k \right) ^2 \, \chi_{J_k \strut}(x) \, dx  \nonumber \\
&=& \sum_{k=1}^{N} \, \int_{J_k} \left(\varphi(x)-
\overline{\varphi}_k \right)^2 \,  dx.  \label{fve1}
\end{eqnarray}
By virtue of Poincar\'{e} inequality we have
\begin{equation}
\int_{J_k} \left(\varphi(x)-
\overline{\varphi}_k \right)^2 \, dx \leq \left(\frac{h}{2 \pi} \right)^2 \int_{J_k} \left(\varphi'(x) \right)^2 \, dx.
\label{poin1}
\end{equation}
Thus, (\ref{fve1}) and (\ref{poin1}) imply
\begin{eqnarray}
\|\varphi(\cdot)-\sum_{k=1}^{N} \, \overline{\varphi}_k \chi_{J_k \strut}(\cdot)\|^2_{L^2}&\leq& \left(\frac{h}{2 \pi} \right)^2 \, \sum_{k=1}^{N} \int_{J_k} \, \left(\varphi'(x) \right)^2 \, dx \nonumber \\
&=&\left(\frac{h}{2 \pi} \right)^2 \int_{0}^{L} \left(\varphi'(x) \right)^2 \, dx,
\end{eqnarray}
which proves inequality (\ref{8}) in the Proposition \ref{prop:89}. \\
Next, we prove inequaliy (\ref{9}). From the Poincar\'{e} inequality (\ref{poin1}) we have
\begin{equation}
\int_{J_k} \, \varphi^2(x) \, dx-
\overline{\varphi}^2_k \, h  \leq \left(\frac{h}{2 \pi} \right)^2 \int_{J_k} \left(\varphi'(x) \right)^2 \, dx.
\end{equation}
Thus, by summing over $k=1,\ldots, N,$ in the above inequality  we conclude inequality
(\ref{9}) of the Proposition \ref{prop:89}.
\begin{thm}
\label{thm:fin-vol}
Let $N$ and $\mu$  be large enough such that $\mu  \geq \nu\left(\frac{2 \pi}{h}\right)^2> \alpha$,
where $\alpha >0$ and $h=\frac{L}{N}$. Then $\|u(t)\|_{L^2}$
tends to zero, as $t \rightarrow \infty$, for every solution $u(t)$
of (\ref{av-RD})-(\ref{av-bc1}).
\end{thm}
{\bf{Proof.}}
We take the $L^{2}$ inner product of equation (\ref{av-RD}) with $u$,
and integrate  by parts to obtain
$$\frac{1}{2} \, \frac{d}{dt} \|u\|_{L^2}^2+\nu \, \|u_x\|_{L^2}^2-\alpha
\|u\|^2_{L^2}+\|u\|^4_{L^4}=-\mu \sum_{j=1}^{N} \frac{L}{N} \, \overline{u}_j^2=
-\mu \frac{L}{N} \,
\gamma^2(u),$$
which implies that
\begin{equation}
\frac{1}{2} \, \frac{d}{dt} \|u\|_{L^2}^2+ \nu \,\|u_x\|_{L^2}^2+\mu \, h\, \gamma^2(u)
-\alpha \, \|u\|^2_{L^2} \leq 0 .
\label{E1}
\end{equation}
Using (\ref{9}), from Proposition \ref{prop:89}, and the assumption that
$\mu  \geq \nu\left(\frac{2 \pi}{h}\right)^2$ we have
\begin{eqnarray}
\nu \,\|u_x\|_{L^2}^2+\mu h  \gamma^2(u)=&\nu \, \left(\frac{h}{2 \pi} \right)^{-2} \left (\left(\frac{h}{2 \pi} \right)^2\|u_x\|_{L^2}^2+h\gamma^2(u)\right) \nonumber \\
&+\left(\mu h \, -\nu\frac{4 \pi^2}{h}\right) \, \gamma^2(u) \nonumber \\
 & \geq \, \frac{4 \pi^2 \, \nu}{h^2} \, \|u\|_{L^2}^2.
\label{E2}
\end{eqnarray}
Substituting (\ref{E2}) in (\ref{E1}) we obtain
$$\frac{1}{2} \, \frac{d}{dt} \|u\|_{L^2}^2+ (\frac{\nu \, 4 \, \pi^2}{h^2}-\alpha) \,\|u\|_{L^2}^2 \leq 0. $$
Therefore, by virtue of Gronwall's inequality and the assumption that
$\nu > \alpha \, \frac{h^2}{4 \, \pi^2}$ one obtains
$$\|u(t)\|_{L^2}^2\, \leq \,  e^{-(\nu (\frac{2 \, \pi \, N}{L})^2-\alpha) \, t} \|u(0)\|_{L^2}^2;$$
and the Theorem follows.

\subsection*{Remark 2.1}
\label{B3}
It is worth mentioning that the assumptions of Theorem \ref{thm:fin-vol}, in particular, that
$N > \sqrt{\frac{L^2 \, \alpha}{4 \pi^2 \nu}}$,    is
consistent  with the fact that the dimension of the unstable manifold about
$\mathbf{v} \equiv 0$ is of order of $\sqrt{\frac{L^2 \, \alpha}{\nu}}$. That is,  one needs at least this number of parameters  to stabilize  $\mathbf{v} \equiv 0$. In Theorem \ref{thm:4.1} we give a different and more general
proof, that
illustrates this point further. As we have mentioned earlier, one can use the same idea to stabilize any other solution, $v(x,t)$, of \eqref{RD}-\eqref{bc1} by using a slightly modified  feedback control in the right-hand side of \eqref{av-RD}-\eqref{av-bc1} of the form $-\mu
\sum_{k=1}^{N} (\overline{u}_k -\overline{v}_k)\; \chi_{J_k\strut} (x)$.


\section{Interpolant operators as feedback controllers}
\label{sec:app-int-fb}
In this section we will consider a general linear map
$I_h:H^1 ([0, L]) \rightarrow L^2([0, L])$ which
is an  interpolant operator that approximates identity with error of order $h$. Specifically, it approximates
the inclusion map $i: H^1 \hookrightarrow L^2$, such  that
the  estimate
\begin{equation}
\|\varphi-I_h (\varphi)\|_{L^2} \leq  \, c \,h \, \|\varphi\|_{H^1},
\label{A1}
\end{equation}
holds, for every $\varphi \in H^1([0, L]).$ The last inequality is a version of the well-known
Bramble-Hilbert inequality, that usually
appears in the context of finite elements \cite{Ciar02}.
We propose here to consider the following general feedback system, to stabilize the solution $v(x,t)$ of \eqref{RD}-\eqref{bc1}, of the form
\begin{eqnarray}
\frac{\partial u}{\partial t} -\nu \, u_{xx}-\alpha \, u+u^3=-\mu \, (I_h(u) -
I_h(v)), \label{rd1:ip} \\
u_x(0)=u_x(L)=0. \label{bc-ip}
\end{eqnarray}
To fix ideas we focus on stabilizing the steady state solution $\mathbf{v} \equiv 0$ of \eqref{RD}-\eqref{bc1}.
Here one can think of $I_h(u)$ as the observables and controllers that will be used
to stabilize our system.\\
Before we state and prove our general theorems concerning system
(\ref{rd1:ip})-(\ref{bc-ip}), we will give some examples of
the approximate  interpolant $I_h(\varphi)$ which satisfy the
approximation property (\ref{A1}).
In particular, we are interested in interpolant operators, $I_h$, of finite-rank, and whose rank is of the order $O(1/h)$.

\subsection{\title Examples of finite-rank approximate identity interpolant operators}
\subsubsection{\title Finite volume elements}
Using the notation of section \ref{sec:fin-vol-feed} we consider the
interpolant operator
\begin{equation}
I_h(\varphi)=\sum_{j=1}^{N} \bar{\varphi}_k \, \chi_{J_k \strut}(x),
\label{ih}
\end{equation}
that uses local spatial averages (finite volume elements) for approximating the local values of the
underlying function. We observe that
the  interpolant operator, $I_h(\varphi)$, that is introduced in  (\ref{ih}) and implemented
in (\ref{rd1:ip}), is exactly the same one discussed in details in section \ref{sec:fin-vol-feed}.
In particular, one can easily see that approximating inequality (\ref{A1}) holds in this
case, thanks to Proposition \ref{prop:89}.

\subsubsection{\title  An interpolant operator based on nodal values}
In this example we consider the  interpolant operator
\begin{equation}
I_h(\varphi)=\sum_{k=1}^{N} \, \varphi(x_k) \, \chi_{J_k \strut} (x), \label{ip-nv}
\end{equation}
where $J_k$ and $\chi_{J_k \strut}$ are as in section \ref{sec:fin-vol-feed}, and the points $x_k \in J_k,$
for $k=1,2,\cdots, N$ are arbitrary. Next, we show that the  interpolant operator given in
(\ref{ip-nv}) satisfied the approximation property (\ref{A1}). Here again we adopt ideas
from \cite{Constantin96,Foias91,Jo:nse2D} to prove the next proposition.
\begin{prop}
\label{prop:nv}
For every $\varphi \in H^1([0, L]) $
$$\|\varphi(\cdot)-\sum_{k=1}^{N} \, \varphi(x_k) \chi_{J_k \strut}(\cdot)\|_{L^2} \, \leq \, h \, \|\varphi_x\|_{L^2} \leq \, h \, \|\varphi\|_{H^1}.$$
\end{prop}
{\bf{Proof.}}
\begin{equation}
\|\varphi(.)-\sum_{k=1}^{N} \, \varphi(x_k) \chi_{J_k \strut}(.)\|^2_{L^2}=\int_{0}^{L} \left(\varphi(x)-\sum_{k=1}^{N}
\, \varphi(x_k) \chi_{J_k \strut} \,(x) \right)^2 \, dx,
\nonumber
\end{equation}
and since $\sum \limits_{k=1}^{N} \chi_{k \strut}(x) \equiv 1$, it follows that
$$
\|\varphi(\cdot)-\sum_{k=1}^{N} \, \varphi(x_k) \, \chi_{J_k \strut}(\cdot)\|^2_{L^2}=\int_{0}^{L} \left(\sum_{k=1}^{N} \, (\varphi(x)-
\, \varphi(x_k)) \, \chi_{J_k \strut} \,(x)) \right)^2 \, dx.
$$
As in the proof of the Proposition \ref{prop:89}, we observe that $\chi_{J_k \strut} (x)\, \chi_{J_l \strut}(x) \equiv \chi_{J_k \strut}(x) \, \delta_{kl}$
and then we obtain
\begin{eqnarray*}
\|\varphi(\cdot)&-&\sum_{k=1}^{N} \, \varphi(x_k) \chi_{J_k \strut}(.)\|^2_{L^2}=\int_{0}^{L} \sum_{k=1}^{N} \left(\varphi(x)-
\, \varphi(x_k) \right)^2 \chi_{J_k \strut} \,(x) \, dx   \\
&=&\sum_{k=1}^{N} \, \int_{J_k} \left(\varphi(x)-
\, \varphi(x_k) \right)^2   \, dx
=\sum_{k=1}^{N} \, \int_{J_k} \left( \int_{x_k}^{x} \varphi'(y) \, dy \right)^2 \, dx   \\
&\leq&\sum_{k=1}^{N} \, \int_{J_k} \left( \int_{J_k} |\varphi'(y)| \, dy \right)^2 \, dx
\leq h \sum_{k=1}^{N} \, \left(\int_{J_k}  |\varphi'(y)| \, dy \right)^2 \, dx.
\end{eqnarray*}
Using Cauchy-Schwarz inequality, we get
\begin{eqnarray*}
\|\varphi(.)-\sum_{k=1}^{N} \, \varphi(x_k) \chi_{J_k \strut}(.)\|^2_{L^2} &\leq&\sum_{k=1}^{N} \, h^2 \, \int_{J_k} \left|\varphi'(y)\right|^2 \, dy, \\
&=& h^2 \|\varphi_x\|_{L^2}^2;
\end{eqnarray*}
which concludes the proof of Proposition \ref{prop:nv}. \\
In view of (\ref{RD}), (\ref{rd1:ip}) and (\ref{ip-nv}) we propose the following feedback
controller for stabilizing $\mathbf{v}\equiv 0$
\begin{eqnarray}
\frac{\partial u}{\partial t}- \nu \, u_{xx}-\alpha u+u^3 =-\mu
\sum_{k=1}^{N} u(x_k)  \;  \chi_{J_k \strut} (x) \label{av-NV} \\
u_x(0)=u_x(L)=0, \label{av-NV-bc1}
\end{eqnarray}
which is a special case of (\ref{rd1:ip}).

\subsubsection{\title Projection onto Fourier modes as an interpolant operator}
Here, we consider the following projection onto the first $N$ Fourier modes as an
example of an  interpolant operator;
\begin{equation}
I_h(\varphi)=\sum_{k=1}^{N} \hat{\varphi}_k \, \cos(\frac{k \, \pi x}{L}), \, \, h=\frac{L}{N},
\label{fm1}
\end{equation}
where the Fourier coefficients are given by
$$\hat{\varphi}_k=\frac{2}{L} \int_{0}^{L} \, \varphi(x) \, \cos(\frac{\pi k x}{L}) dx. $$
Next, we observe that inequality (\ref{A1}) holds for the  interpolant operator given in (\ref{fm1}). \\
\begin{prop}
\label{prop:fm}
Let $\varphi \, \in H^1([-L, L])$ be an even function, i.e. $\varphi(-x)=\varphi(x).$
Then
\begin{equation}
\| \varphi(x)-\sum_{k=1}^{N} \, \widehat{\varphi}_k \,
\cos \left(\frac{k \, x  \pi} {L} \right)\|_{L^2([0, L])} \leq c \, h \, \|\varphi_x\|_{L^2([0, L])}.
\label{fm2}
\end{equation}
\end{prop}
{\bf{Proof.}}
The proof of this proposition is a simple exercise in Fourier series.
Thus, it will be omitted.\\


\section{Existence, uniqueness and  stabilization using the $I_h$ feedback control}
\label{sec:exis-uni}
In this section we establish the global existence and uniqueness for the general feedback system
introduced in (\ref{rd1:ip})-(\ref{bc-ip}); and that the $I_h$ feedback control is stabilizing the steady state solution $v\equiv 0$ of the \eqref{RD}-\eqref{bc1}. This will be accomplished under the assumptions
(\ref{A1}) and that $\mu$ is large enough, and $h$ is small enough,   satisfying:
\begin{equation}
\mu \ge (2\alpha + 3\nu L^{-2}) \qquad \hbox{and} \quad
\nu \geq \mu \, c^2 \, h^2.
\label{assump.36}
\end{equation}
To this end one uses the standard Galerkin approximation
procedure based on the eigenfunctions of the Laplacian,
subject to the Neumann boundary
condition, i.e., ${\cos(\frac{\pi \, k\, x}{L})}$ for $k=1,2...$. We
will omit the details
of this standard procedure and provide only the formal \textit{a-priori} estimates
(see, e.g., \cite{Tem97}). These estimates
can be  obtained rigorously through the Galerkin procedure, by passing to the limit
while using the relevant compactness theorems. \\
Let us now establish the aformentioned formal {\it a-priori} bounds for the solution which are
essential for guaranteeing globlal existence and uniqueness. \\
System (\ref{rd1:ip})-(\ref{bc-ip}) can be rewritten as
\begin{eqnarray}
\frac{\partial u} {\partial t} - \nu \, u_{xx}+\frac{\nu}{L^2} u -(\alpha+\frac{\nu}{L^2}) \, u=-u^3-\mu I_h(u) \label{geu1} \\
u_x(0)=u_x(L)=0.
 \label{geu2}
\end{eqnarray}
Taking the $L^2$- inner product of (\ref{geu1}) with $u$, integrating by
parts and using the Neumann boundary
conditions, we obtain
\begin{eqnarray}
\frac{1}{2} \frac{d} {d t}  \int_{0}^{L} u^2 \, dx &+& \nu \,
\int_{0}^{L} u^2_{x} \, dx + \frac{\nu}{L^2} \int_{0}^{L} u^2
\, dx= -\int_{0}^{L} u^4 dx \nonumber \\
& & +(\alpha+\frac{\nu}{L^2}) \int_{0}^{L} u^2 \, dx-\mu \,
\int_{0}^{L} I_h(u) \, u \, dx. \nonumber
\end{eqnarray}
Writing
$$
I_h(u) \, u= \left(I_h(u)-u \right)u + u^2
$$
and applying the Cauchy-Schwarz  inequality, we get
\begin{equation}
\frac{1}{2} \frac{d} {dt}  \int_{0}^{L} u^2 \, dx+ \nu \,
\int_{0}^{L} u^2_{x} \, dx + \frac{\nu}{L^2} \int_{0}^{L} u^2
\, dx \leq -   \int_{0}^{L} u^4 dx+
(\alpha+\frac{\nu} {L^2}) \|u\|^2_{L^2}
\nonumber
\end{equation}
\begin{equation}
-\mu \, \int_{0}^{L} u^2 \, dx+\mu \,
\left(\int_{0}^{L} u^2 \, dx\right)^\frac{1}{2}
\, \left(\int_{0}^{L} |u-I_h(u)|^2 \, dx\right)^\frac{1}{2}.
 \nonumber
\end{equation}
Using  Young's inequality  we reach
\begin{eqnarray}
\frac{1}{2} \frac{d} {dt}  \int_{0}^{L} u^2 \, dx+ \nu \,
\int_{0}^{L} u^2_{x} \, dx &+& \frac{\nu}{L^2}
\int_{0}^{L} u^2 \, dx \leq -
\int_{0}^{L} u^4 dx+(\alpha+\frac{\nu} {L^2}) \|u\|^2_{L^2} \nonumber \\
&-&\frac{\mu}{2} \, \int_{0}^{L} |u|^2 \, dx+
\frac{\mu}{2} \|u-I_h(u)\|^2_{L^2} \, dx.
 \nonumber
\end{eqnarray}
Using (\ref{A1}), and the definition of the $H^1$-norm given
in (\ref{norm-h1})
we obtain
\begin{equation}
\frac{1}{2} \frac{d} {d t}  \int_{0}^{L} u^2 \, dx+ \nu \,
\int_{0}^{L} u^2_{x} \, dx + \frac{\nu}{L^2} \int_{0}^{L} u^2
\, dx \leq -  \int_{0}^{L} u^4 dx+
(\alpha+\frac{\nu} {L^2}) \|u\|^2_{L^2} \nonumber
\end{equation}
\begin{equation}
-\frac{\mu}{2} \, \int_{0}^{L} |u|^2 \, dx +\mu \;
\frac{c^2 \, h^2}{2} \left(\frac{1}{L^2} \int_{0}^{L} u^2 \, dx
+ \int_{0}^{L} u^2_{x} \, dx \right). \nonumber
\end{equation}
Thanks to the assumption (\ref{assump.36}) we conclude
\begin{equation}
\frac{d}{dt} \, \|u\|_{L^2}^2+\nu\, (\|u_x\|_{L^2}^2+\frac{1}{L^2} \, \|u\|_{L^2}^2)
 \; \leq 0.   \label{geu3}
\end{equation}
Therefore, by dropping the $\|u_x\|_{L^2}^2$ term from the left-hand side of (\ref{geu3})
and applying Gronwall's inequality we have
\begin{equation}
\|u(t) \|^2_{L^2} \leq e^{\frac{-\nu \, t}{L^2} } \, \|u(0)\|^2_{L^2}=:K_0(t).
\label{E8}
\end{equation}
Notice that from \eqref{geu3} one also concludes that for every $\tau>0$
\begin{equation}
\nu\, \int_t^{t+\tau}(\|u_x(s)\|_{L^2}^2+\frac{1}{L^2} \, \|u(s)\|_{L^2}^2)\, ds
 \; \leq \|u(t)\|^2_{L^2}.   \label{H1Es}
\end{equation}
Next, we show the continuous dependence of the solutions of (\ref{rd1:ip}) on the initial data and the
uniqueness, provided the assumptions (\ref{A1}) and (\ref{assump.36}) hold.
Indeed, let $u_1, u_2$ be two solutions and
$w=u_1-u_2$ of (\ref{rd1:ip}). From (\ref{rd1:ip}) we find that
$$\frac{\partial w}{\partial t} -\nu w_{xx}-\alpha\, w=u^3_2-u^3_1-\mu\, I_h(w).$$
Multiplying by $w$ and integrating  with respect to $x$ over $[0, L]$ we get
\begin{eqnarray}
\frac{1}{2} \frac{d} {d t}  \int_{0}^{L} w^2 \, dx + \nu \, \int_{0}^{L} w^2_{x} \, dx&=&\alpha \,
 \int_{0}^{L} w^2 \, dx - \int_{0}^{L} \, w^2 \, \frac{(u_1+u_2)^2 +u_1^2+u_2^2}{2} \, dx \nonumber \\
 &&-\mu \,  \int_{0}^{L} I_h(w) \, w \ dx \nonumber  \\
&\leq& (\alpha-\mu) \, \int_{0}^{L} w^2 \, dx+\mu \, \int_{0}^{L} |I_h(w)-w| \, |w| \, dx. \nonumber
\end{eqnarray}
A straightforward computation, using the Cauchy-Schwarz and Young inequalities and
assumption (\ref{A1}),  yields
\begin{eqnarray}
\frac{1}{2} \frac{d} {d t}  \int_{0}^{L} w^2 \, dx + \nu \, \int_{0}^{L} w^2_{x} \, dx&\leq&
(\alpha-\mu) \, \int_{0}^{L} w^2 \, dx+\mu \, \|I_h(w)-w\|_{L^2} \, \|w\|_{L^2}  \nonumber \\
&\leq& (\alpha-\mu) \, \|w\|^2_{L^2} \, dx+\frac{\mu}{2} \, \|I_h(w)-w\|^2_{L^2}+\frac{\mu}{2} \|w\|^2_{L^2}  \nonumber \\
&\leq& (\alpha-\frac{\mu}{2}) \, \|w\|^2_{L^2}+\frac{\mu}{2} \, c^2 \, h^2 \|w\|^2_{H^1}.  \nonumber
\end{eqnarray}
Using (\ref{norm-h1}), the definition of the $H^1$-norm, we reach
\begin{eqnarray}
\frac{1}{2} \frac{d} {d t}  \int_{0}^{L} w^2 \, dx + \nu \, \int_{0}^{L} w^2_{x} \, dx
&&\leq \left(\alpha-\frac{\mu}{2}\right) \, \|w\|^2_{L^2}+\frac{\mu}{2} \, c^2 \, h^2 \,
\left(\frac{\|w\|^2_{L^2}}{L^2}+\|w_x\|^2_{L^2} \right)  \nonumber \\
 &\leq& \left(\alpha-\frac{\mu}{2}+\frac{\mu}{2} \, \frac{c^2 \, h^2}{L^2}\right) \, \|w\|^2_{L^2}+\frac{\mu}{2} \,
c^2 \, h^2 \, \|w_x\|^2_{L^2}.  \nonumber
\end{eqnarray}
By assumption (\ref{assump.36}) the above implies
$$
\frac{1}{2} \frac{d} {d t}  \|w \|^2_{L^2}+ \frac{\nu}{2} \, \|w_{x}\|^{2}_{L^2} \leq \left(\alpha-\frac{\mu}{2}
+\frac{\nu} {2 L^2} \right) \, \|w\|^{2}_{L^2}.
$$
Using assumption (\ref{assump.36}) one more time the above inequality simplifies to
$$
\frac{1}{2} \frac{d} {d t}  \|w \|^2_{L^2}+ \frac{\nu}{2} \, \|w_{x}\|_{L^2}^{2} \leq - \frac{\nu} { L^2} \|w\|_{L^2}^2.
$$
Finally, by Gronwall's inequality we have
\begin{equation}
\|w(t)\|_{L^2}^2 \leq e^{- \frac{\nu t} { L^2}} \, \|w(0)\|_{L^2}^2.
\label{E10}
\end{equation}
Thus, if $w(0)=0$ then $\|w(t)\|_{L^2} \equiv 0$. Moreover, inequaliy (\ref{E10}) implies
the continuous dependence of the solutions of (\ref{rd1:ip})-(\ref{bc-ip}) on the
initial data. In conclusion, from the above, and in particular thanks to (\ref{E8})
and (\ref{H1Es}), we have the following theorem:
\begin{thm}
\label{thm:4.1}
Let $\mu, \nu$ and $h$ be positive parameters satisfying assumption (\ref{assump.36});
and that
$I_h$ satisfies (\ref{A1}). Suppose $T \, > \, 0$ and $u_0 \in L^2([0, L]),$ then system
(\ref{rd1:ip})-(\ref{bc-ip})  has a unique solution $u \in C([0, T], L^2) \cap L^2([0, T], H^1)$
which also depends continuously on the initial data. Moreover,
$$\lim_{t \rightarrow \infty} \|u (t)\|_{L^2}^2 =0;
$$
and for every $\tau >0$
$$\lim_{t \rightarrow \infty} \int_t^{t+\tau}\|u_x(s)\|_{L^2}^2 \, ds =0.$$
In particular, we concluded the under that above assumption the feedback control interpolant operator $I_h$ is stabilizing the steady state solution
$v\equiv 0$ of \eqref{RD}-\eqref{bc1}.
\end{thm}

\subsection*{Remark 4.1}
Let us observe that in order to satisfy assumption (\ref{assump.36}) one can choose, for small values of $\nu$,
$\mu=O(\alpha)$. As a result, assumption (\ref{assump.36}) will hold if we
choose $h$ small enough such that $N:=\frac{L}{h}=O(\sqrt{\frac{\alpha \, L^2}{\nu}})$, that is the number of
feedback controllers  is comparable to the dimension of the unstable manifold
about $\mathbf{v} \equiv0$. This is consistent with our earlier observation in
the introduction and in Remark 2.1.

\section{Stabilizing in the $H^1$-norm}
\label{sec:stabinh1}
In the previous section we have shown that the feedback system (\ref{rd1:ip})-(\ref{bc-ip})
stabilzes the steady state solution $\mathbf{v} \equiv 0$ in the $L^2$-norm, i.e.,
$\|u\|_{L^2} \rightarrow 0$, as $t \rightarrow \infty$, provided assumptions
(\ref{assump.36})  holds.
Next, we show that we also have $\|u(t)\|_{H^1} \rightarrow 0$, as $t \rightarrow \infty$. To this
end it is enough to show that $\|u_x\|_{L^2} \rightarrow 0,$ as $t \rightarrow \infty$. \\
Let us rewrite (\ref{rd1:ip})-(\ref{bc-ip}) as
\begin{eqnarray}
u_t+\frac{\nu}{L^2} \, u-\nu \, u_{xx}-(\alpha+\frac{\nu}{L^2}) \, u+u^3=-\mu \, I_h(u) \label{rd:ip-st-2} \\
u_x(0)=u_x(L)=0.
 \label{rd:ip-st-bc2}
\end{eqnarray}
Thanks to the estimate \eqref{H1Es} we realize that the solution instantaneously becomes in $H^1$. Therefore, without loss of generality we can assume that the initial data $u_0 \in H^1$. Below we provide the formal arguments and  estimates, which, as we have already indicated earlier, can be established rigorously by  using a Galerkin approximation procedure.  We take the $L^2$ inner product of (\ref{rd:ip-st-2}) with $-u_{xx}$, integrating by parts, and using the
Neumann boundary conditions  (\ref{bc-ip}) we obtain:
\begin{eqnarray*}
\frac{1}{2} \frac{d} {d t}  \|u_x \|^2_{L^2}&+&\nu \, \|u_{xx}\|_{L^2}^2+ \frac{\nu}{L^2} \, \|u_x\|^2_{L^2}
-(\alpha+\frac{\nu}{L^2}) \|u_x\|^2   \\
&=& \int_{0}^{L} u^3 \, u_{xx} \, dx  + \mu \, \int_{0}^{L} I_h(u) \, u_{xx} \, dx    \\
&=&-3 \, \int_{0}^{L} u^2 \, u_{x}^2 \, dx+ \mu \, \int_{0}^{L} (I_h(u)-u) \, u_{xx} \, dx
 + \mu \, \int_{0}^{L} u \, u_{xx} \, dx \\
&=&-3 \, \int_{0}^{L} u^2 \, u_{x}^2 \, dx+ \mu \, \int_{0}^{L} (I_h(u)-u) \, u_{xx} \, dx
 - \mu \, \int_{0}^{L} \, u_{x}^2 \, dx.
\end{eqnarray*}
By Cauchy-Schwarz inequality we have
\begin{eqnarray}
\frac{1}{2} \frac{d} {d t}  \|u_x \|^2_{L^2}+ \nu \, \|u_{xx} \|_{L^2}^2 +\frac{\nu}{L^2} \|u_x\|^2_{L^2}
 &\leq& (\alpha+\frac{\nu}{L^2})  \, \|u_x\|^2_{L^2}
-\mu \, \|u_x\|^2_{L^2} \nonumber \\ &+&\mu \, \|I_h(u)-u)\|_{L^2} \, \|u_{xx}\|_{L^2}. \nonumber
\end{eqnarray}
Applying Young's inequality we obtain
\begin{equation}
\frac{1}{2} \frac{d} {d t}  \|u_x \|^2_{L^2}+ \nu \, \|u_{xx} \|^2_{L^2}
 +\frac{\nu}{L^2} \|u_x\|^2_{L^2} \leq (\alpha+\frac{\nu}{L^2}-\mu) \|u_x\|^2
+\frac{\nu}{2}\, \|u_{xx}\|^2_{L^2}+\frac{\mu^2}{2\, \nu} \, \|I_h(u)-u)\|_{L^2}^2. \nonumber
 \nonumber \\
\end{equation}
Using property (\ref{A1}) and the definition of the $H^1-$norm in
(\ref{norm-h1}), we have
\begin{eqnarray}
\frac{1}{2} \frac{d} {d t}  \|u_x \|^2_{L^2} \, dx + \frac{\nu}{2} \, \|u_{xx} \|^2_{L^2} + \frac{\nu}{L^2} \|u_x\|_{L^2}^2 &\leq& (\alpha+\frac{\nu}{L^2}-\mu) \|u_x\|^2_{L^2} \nonumber \\
&+& \frac{h^2 \, \mu^2 \, c^2}{2\, \nu} \, \left( \|u_x\|_{L^2}^2+\frac{1}{L^2} \|u\|_{L^2}^2 \right). \nonumber
\end{eqnarray}

\begin{eqnarray}
\frac{1}{2} \frac{d} {d t}  \|u_x \|^2_{L^2} \, dx + \frac{\nu}{2} \, \|u_{xx} \|^2_{L^2} + \frac{\nu}{L^2} \|u_x\|_{L^2}^2 &\leq& (\alpha+\frac{\nu}{L^2} +\frac{h^2 \, \mu^2 \, c^2}{2\, \nu}-\mu) \|u_x\|^2_{L^2} \nonumber \\
&+& \frac{h^2 \, \mu^2 \, c^2}{2\, \nu \, L^2} \, \|u\|_{L^2}^2. \nonumber
\end{eqnarray}

Thanks to assumption (\ref{assump.36}) we observe that
$$
\frac{h^2 \, \mu^2 \, c^2}{2\, \nu} \le \frac{\mu}{2} \quad \hbox{and} \quad (\alpha+\frac{\nu}{L^2} +\frac{h^2 \, \mu^2 \, c^2}{2\, \nu}-\mu) \le -\frac{\nu}{2\, L^2}.
$$

Therefore, the above implies

\begin{equation}
\frac{1}{2} \frac{d} {d t}  \|u_x \|^2_{L^2}+ \frac{3 \, \nu}{2 \, L^2} \|u_x\|_{L^2}^2
\leq \frac{\mu}{2} \, \|u\|_{L^2}^2. \nonumber
\end{equation}
Since $\lim_{t \rightarrow \infty} \|u(t)\|^2=0,$
then by Gronwall's inequality it is easy to show that $\|u_x\|_{L^2}^2 \rightarrow 0$, as $t \rightarrow \infty$,
 (see also special Gronwall's type Lemma in \cite{Jo:nse2D}).

\section{Nodal observables  and feedback controllers}
In this section we propose a different feedback control based on nodal value observables and feedback
controllers. Assume that the observables are the values of the solutions $u(\overline{x}_k)$, at the points
$\overline{x}_k \, \in \, J_k=[(k-1)\, \frac{L}{N}, k \, \frac{L}{N}], \, k=1,...,N,$ and that the feedback
is at some points $x_k\, \in J_k, \, x_k$ is not necessarily the same as $\overline{x}_k$. That is the
measurements are made at $\overline{x}_k$, while the feedback controllers are at $x_k$, for
$k=1,2...,N.$ To avoid technical issues that are dealing with boundary conditions, we focus here
on the periodic boundary condition case. In this case the feedback system will read
\begin{equation}
\frac{\partial u}{ \partial t}-\nu \, u_{xx}-\alpha u+u^3=-\mu \; \sum_{k=1}^{N}
\, h \, u(\overline{x}_k) \, \delta(x-x_k),
\label{71}
\end{equation}
\begin{equation}
u(x, t)=u(x+L, t),
 \label{72}
\end{equation}
where $h=\frac{L}{N}$; and $\delta(x-a) \, \in \, H^{-1}_{per}([0, L])$,  for $a \, \in \, [0, L]$,
and is extended periodically such that
\begin{equation}
<\delta(\cdot-a), \varphi>=\varphi(a)
\label{73}
\end{equation}
for every $\varphi \, \in H^{1}_{per}([0, L])$. \\
The feedback control proposed in (\ref{71}) -(\ref{72}) is different than that of (\ref{rd1:ip})-(\ref{bc-ip}),
since the right-hand side in (\ref{71}) is a distribution that belongs to $H^{-1}_{per}([0, L])$, while
the right-hand side in (\ref{rd1:ip}) belongs to $L^2([0, L]).$ \\
In this section we will show that, under similar assumptions to those in Theorem \ref{thm:4.1},
the proposed feedback system (\ref{71}) stabilizes the steady state $\mathbf{v} \equiv 0$ in the
$L^2-$norm. One should not expect here a stronger statement, as the one stated in section
\ref{sec:stabinh1}, in which the stabilizing is also valid in the $H^1-$norm. This is because
the solutions of (\ref{71})-(\ref{72}) are weaker
than those of (\ref{rd1:ip})-(\ref{bc-ip}), since the right-hand side
in (\ref{71}) is less regular than its counterpart in (\ref{rd1:ip}).\\
 As before, we will show, below,  the formal steps, which demonstrate simultaneously the global existence, uniqueness and stabilizing effect.
 These formal steps and estimates can be justified rigorously by implementing
the Galerkin procedure based on the eigenfunction of the Laplacian, subject to periodic
boundary conditions, with period $L$ (see, e.g., \cite{Tem97}).
First, let us prove the following Lemma, which is basically the embedding of the H\"{o}lder space
of $C^{\frac{1}{2}} \, \subset H^1$ (see, e.g.,  \cite{Constantin96,Tem97}).

\begin{lem}
\label{lem-sec71}
Let $x_k, \overline{x}_k  \in J_k=[(k-1) \, h, k \, h], k=1,..,N, $
where $h=\frac{L}{N}$, $N \in \NN.$
Then for every $\varphi \, \in \, H^1([0, L])$ we have
\begin{equation}
\sum_{k=1}^{N} | \varphi(x_k)-
\varphi(\overline{x}_k)|^ 2 \leq \, h \, \|\varphi_x\|^2_{L^2},
\label{7lem1}
\end{equation}
and
\begin{equation}
\|\varphi\|^2_{L^2}  \leq \, 2 \, \left[h \,  \sum_{k=1}^{N} |\varphi(x_k)|^2+
h^2 \, \|\varphi_x\|^2_{L^2} \right].
\label{7lem2}
\end{equation}
\end{lem}
{\bf{Proof.}}
We prove inequality (\ref{7lem1}) for $\varphi \, \in \, C^1([0, L])$, and by the density
of $C^1 \subset H^1$
the result follows for every $\varphi \, \in \,H^1$.
\begin{eqnarray}
\left| \varphi(x_k)- \varphi(\overline{x}_k)\right|^ 2 &\leq&
\left|\int_{\overline{x}_k}^{x_k}\varphi'(s) \, ds\right|^2 \,
\leq \, \left(\int_{J_k} |\varphi'(s)| \, ds\, \right)^2 \nonumber \\
&\leq& |J_k| \, \int_{J_k} |\varphi'(s)|^2 \, ds=h \, \int_{J_k} |\varphi'(s)|^2 \, ds.
\nonumber
\end{eqnarray}
By summing the above inequality over $k=1,..,N$ we conclude (\ref{7lem1}).\\
To prove (\ref{7lem2}) we observe that for every $x \, \in \, J_k$ we have
\begin{equation}
|\varphi(x)| \, \leq \, |\varphi(x_k)|+\int_{J_k} |\varphi'(s)| \, ds.
\nonumber
\end{equation}
Thus
\begin{equation}
|\varphi(x)|^2 \, \leq \, 2\,  \left[|\varphi(x_k)|^2
+\left(\int_{J_k} |\varphi'(s)| \, ds \right)^2 \right],
\end{equation}
and by integrating with respect to $x$ over $J_k$, and using the Cauchy-Schwarz inequality, we obtain
\begin{equation}
\int_{J_k} \, |\varphi(x)|^2 \, dx\leq \, 2\, h \,\left[|\varphi(x_k)|^2+h \, \int_{J_k} \, |\varphi'(s)|^2 \, ds \,\right].
\end{equation}
Now we conclude (\ref{7lem2}) by summing over $k=1,\dots, N.$

\begin{thm}
\label{thm}
Let $\mu\, > \, 4 \, \alpha$ and $h$ is small enough such that $\nu \geq 2 \, \mu \, h^2.$
Then for every $T>0,$ and every $u_0 \in L^2_{per}[0, L]$ system (\ref{71}) has a unique solution
$$u \, \in \, C([0, T]; L^2_{per}[0, L]) \cap  L^2 \left([0, T]; H^1_{per}\right[0, L]) \cap
L^4\left([0, T]; L^4_{per}[0, L]\right),$$
and
$$\frac{\partial u} {\partial t} \, \in \, L^2\left([0, T]; H^{-1}_{per} \right).$$
Moreover,
\begin{equation}
\lim_{t \rightarrow \infty} \|u(t)\|_{L^2}=0, \quad \hbox{and for every} \quad \tau >0 \quad  \lim_{t \rightarrow \infty} \int_t^{t+\tau} \|u_x(s)\|^2_{L^2} \, ds=0
\label{74}
\end{equation}
\end{thm}
{\bf{Proof.}}
We will show here only the relevant {\it a priori} estimates. The rest of the regularity results are standard for nonlinear parabolic equations (see, e.g., \cite{Tem97}). \\
We take the $H^{-1}$ action of (\ref{71}) on $u \in H^{1}$, and use Lemma of Lions-Magenes
 (cf. Chap. III-p.169, \cite{Tem01}),
to obtain:
\begin{eqnarray}
\frac{1}{2} \frac{d}{dt} \|u\|^2_{L^2} &+&\nu \, \|u_x\|_{L^2}^2=\alpha \, \| u\|^{2}_{L^2}-\int_{0}^{L} \, u^4 \, dx
-\mu \, h \, \sum_{k=1}^{N} \, u(\overline{x}_k) \, u(x_k)  \nonumber\\
&=&\alpha \, \| u\|^{2}_{L^2}-\int_{0}^{L} \, u^4 \, dx
-\mu \, h \, \sum_{k=1}^{N} \, |u(x_k)|^2
+ \mu \, h \, \sum_{k=1}^{N} \, \left(u(x_k)-u(\overline{x}_k) \right) \, u(x_k) \nonumber \\
&\leq&\alpha \, \| u\|^{2}_{L^2}-\int_{0}^{L} \, u^4 \, dx
-\frac{\mu}{2} \, h \, \sum_{k=1}^{N} \, |u(x_k)|^2
+\frac{\mu}{2} \, h \, \sum_{k=1}^{N} \, \left|u(x_k)-u(\overline{x}_k) \right|^2, \nonumber
\end{eqnarray}
where in the last step we applied the Young's inequality. Next, we apply (\ref{7lem1}) and (\ref{7lem2})
to the right-hand side
$$
\frac{1}{2} \frac{d}{dt} \|u\|^2_{L^2}+\nu \, \|u_x\|_{L^2}^2 \leq \,
\alpha \, \| u\|^{2}_{L^2}-\int_{0}^{L} \, u^4 \, dx
-\frac{\mu}{4}  \, \|u\|^2_{L^2}+\mu \, h^2 \, \|u_x\|_{L^2}^{2}.
$$
Hence
$$
\frac{1}{2} \frac{d}{dt} \|u\|^2_{L^2} +(\nu-\mu \, h^2) \, \|u_x\|_{L^2}^2
\leq \, (\alpha-\frac{\mu}{4})\, \| u\|^{2}_{L^2}-\int_{0}^{L} \, u^4 \, dx.
$$
Since $\nu \geq 2\mu \, h^2$ and $4 \alpha < \mu$ we conclude:
$$
 \frac{d}{dt} \|u\|^2_{L^2} +\nu\, \|u_x\|_{L^2}^2 + (\frac{\mu}{2}- 2 \alpha)\, \| u\|^{2}_{L^2}
\leq 0.
$$
Thanks
to Gronwall's inequality  we conclude  from the  above \eqref{74}, and the regularity of the solutions as stated in the theorem.

 Next, we prove the uniqueness of solutions.
Let $u_1$ and $u_2$ be any two solutions. Denote by $w=u_1-u_2$. Then $w$
satisfies
\begin{equation*}
\frac{\partial w}{\partial t}-\nu \, w_{xx}-\alpha \, w+(u_1^2+u_1 \, u_2+u_2^2) \, w=-\mu \, h \,
 \sum_{k=1}^{N} \, w(\overline{x}_k) \, \delta(x-x_k).
\end{equation*}
Taking the $H^{-1}$ action on $w \in \, H^1$, and using again Lemma of Lions-Magenes
(cf. Chap. III-p.169, \cite{Tem01}), we obtain
\begin{equation}
\frac{1}{2} \frac{d}{dt} \|w\|^2_{L^2} +\nu \, \|w_x\|_{L^2}^2
-\alpha \, \| w\|^{2}_{L^2}= - \int_{0}^{L} \, (u_1^2+u_1 \, u_2+u_2^2) \, w^2 \, dx
-\mu \, h \, \sum_{k=1}^{N} \, w(\overline{x}_k) \, w(x_k). \nonumber
\end{equation}
Since $\int_{0}^{L} \, (u_1^2+u_1\, u_2+u_2^2) \, w^2 \, dx = \int_{0}^{L} \,
\left (\frac{u_{1}^{2}+u_2^2}{2} + \left( \frac{u_1+u_2}{2}  \right)^2 \right )\, w^2 \, dx \geq 0$ we obtain
\begin{equation}
\frac{1}{2} \frac{d}{dt} \|w\|^2_{L^2} +\nu \, \|w_x\|_{L^2}^2
-\alpha \, \| w\|^{2}_{L^2} \leq -\mu \, h \, \sum_{k=1}^{N} \, w(\overline{x}_k) \, w(x_k). \nonumber
\end{equation}
Next,  we follow the same steps as in the beginning of the proof   to obtain
\begin{equation}
\frac{1}{2} \frac{d}{dt} \|w\|^2_{L^2} +(\nu-\mu \, h^2) \, \|w_x\|_{L^2}^2
\leq (\alpha-\frac{\mu}{4}) \, \| w\|^{2}_{L^2}.  \nonumber
\end{equation}
Since $\nu \geq 2\mu \, h^2$ and $\mu \, > \, 4 \, \alpha$ we conclude, thanks
to Gronwall's inequality,
\begin{equation}
\|w(t)\|^2_{L^2} \leq e^{(\alpha-\mu/4) \, t} \|w(0)\|^2_{L^2}.
\label{76}
\end{equation}
Notice that (\ref{76}) implies the uniqueness of the solutions and their continuous
 dependence on the initia data.
\subsection*{Remark 6.1}
Here again we observe that for small values of $\nu$ by choosing $\mu=O(\alpha)$ then the condition of the
theorem imply that $N:=\frac{L}{h}=O(\sqrt{\frac{\alpha}{\nu}})$; which is comparable
to the dimension of the unstable manifold about the steady state $\mathbf{v} \equiv 0.$


 \subsection*{Acknowledgements}
 The work of A.A.~is supported in part by the DFG grants SFB-910 and SFB-947.
The work of E.S.T.~is  supported in part by the   NSF  grants DMS-1009950, DMS-1109640,  and DMS-1109645, and by the Alexander von Humboldt Stiftung/Foundation and the
Minerva Stiftung/Foundation.  Partial support was also provided by the  CNPq-CsF grant \# 401615/2012-0, through the program Ci\^encia sem Fronteiras.
 E.S.T. is also thankful to the kind hospitality  of the Freie
Universit\"at Berlin, where this work was initiated, and of the Instituto Nacional de Matem\' {a}tica  Pura  e Aplicada (IMPA), Brazil,  where part of this work was completed.

\end{document}